\documentclass{article}
\usepackage{amssymb}

\newcommand{\Tilde}{\widetilde}
\newcommand{\R}{{\Bbb R}}

\newcommand{\C}{{\Bbb C}}

\newcommand{\Z}{{\Bbb Z}}
\newcommand{\Q}{{\Bbb Q}}

\newcommand{\CP}{{\Bbb CP}}
\newcommand{\Dd}{{\mathcal D}}
\newcommand{\Aa}{{\mathcal A}}
\newcommand{\Ff}{{\mathcal F}}
\newcommand{\Mm}{{\mathcal M}}
\newcommand{\Nn}{{\mathcal N}}

\newcommand{\Jj}{{\mathcal J}}
\newcommand{\Ss}{{\mathcal S}}
\newcommand{\Xx}{{\mathcal X}}

\newcommand{\SO}{{\rm SO}}
\newcommand{\Aut}{{\rm Aut}}
\newcommand{\End}{{\rm End}}
\newcommand{\Gr}{{\rm Gr}}

\newcommand{\p}{{\partial}}
\newcommand{\al}{{\alpha}}

\newcommand{\be}{{\beta}}

\newcommand{\Om}{{\Omega}}
\newcommand{\om}{{\omega}}
\newcommand{\eps}{{\varepsilon}}
\newcommand{\de}{{\delta}}

\newcommand{\ga}{{\gamma}}

\newcommand{\io}{{\iota}}
\newcommand{\ka}{{\kappa}}
\newcommand{\la}{{\lambda}}
\newcommand{\La}{{\Lambda}}
\newcommand{\si}{{\sigma}}
\newcommand{\Si}{{\Sigma}}
\newcommand{\Symp}{{\rm Symp}}
\newcommand{\Diff}{{\rm Diff}}

\newcommand{\Hom}{{\rm Hom}}
\newcommand{\Map}{{\rm Map}}

\newcommand{\Fol}{{\rm Fol}}

\newcommand{\SmS}{{\smallskip}}
\newcommand{\MS}{{\medskip}}

\newcommand{\NI}{{\noindent}}
\newcommand{\proof}[1]{\noindent{\bf Proof#1:\  }}
\newcommand{\QED}{\hfill$\Box$\medskip}
\newtheorem{theorem}{Theorem}[section]
\newtheorem{thm}[theorem]{Theorem}
\newtheorem{cor}[theorem]{Corollary}
\newtheorem{defn}[theorem]{Definition}
\newtheorem{rmk}[theorem]{Remark}
\newtheorem{lemma}[theorem]{Lemma}

\newtheorem{prop}[theorem]{Proposition}

\title{Symplectomorphism Groups and Almost Complex Structures}
\author{
Dusa McDuff\thanks{Partially
supported by NSF grants DMS 9704825 and 0072512. Keywords: symplectomorphism group,
symplectic ruled surface, J-holomorphic curve.  
MSC 2000: Primary 57R17, 57S05; Secondary 53D35} \\ State
University of New York at Stony Brook, USA \\ (dusa@math.sunysb.edu)}

\date{August 29, 2000, last revised December 26, 2004}

\begin{document}

\maketitle
\MS

\begin{center}for Andr\'e Haefliger's 70th
birthday
\end{center}
\MS\MS

\MS

\NI
{\bf Abstract}
This paper studies groups  of
symplectomorphisms of  ruled surfaces $M$ for 
symplectic forms with varying cohomology class.  This cohomology class is
characterised by the ratio $\mu$ of the size of the base to that of the fiber. 
By considering appropriate spaces of almost complex structures, we
investigate how the topological type of these groups  changes as $\mu$
increases.  If the base is a sphere, this changes precisely when $\mu$ passes
an integer, and, for general bases, it stabilizes as $\mu \to \infty$.   Our
results extend and make more precise some of the conclusions of
Abreu--McDuff concerning the rational homotopy type of these groups for
rational ruled surfaces. \MS\MS

\section{Introduction}

One of the interesting facts of symplectic geometry is that symplectic
manifolds admit a large family of symplectomorphisms, i.e.
diffeomorphisms that preserve the symplectic structure.  Indeed, 
 every function $H$  on a
closed symplectic manifold $(M, \om)$ generates a flow
$\phi_t^H, t\in \R,$ consisting of   symplectomorphisms. 
Thus the group $\Symp(M,\om)$ of 
all  symplectomorphisms of $(M, \om)$ is always infinite dimensional.
On the other hand,  in some special cases the homotopy
type of this group can be calculated and  it turns out to be not too large.
For example, although nothing at all is known about the group of compactly
supported diffeomorphisms of $\R^4$ --- it is not even known whether it is
connected --- Gromov showed in 1985 (see~\cite{G}) that the group of
compactly supported symplectomorphisms of $\R^4$ with its standard
symplectic structure is contractible.  He also showed that the group of
symplectomorphisms of $S^2\times S^2$, with the product symplectic
form in which both spheres have the same size, has the homotopy type of
a Lie group:  in fact it is homotopic to the semi-direct product of
$\SO(3)\times \SO(3)$ with $\Z/2\Z$.

More recently, Abreu~\cite{A}  showed that
when one sphere factor is larger than the other the group of
symplectomorphisms is no longer homotopy equivalent to a Lie group, 
because it does not have the right kind of rational homotopy type.  
This rational homotopy type was calculated in 
Abreu--McDuff~\cite{AM}, and was found to change precisely when the ratio
$\mu$ of the size of the larger to the smaller sphere passes an integer value.  

The current paper extends this last result to the
actual rather than the rational homotopy type.  Our 
arguments are more general than before in that in principle they apply to
ruled surfaces over Riemann surfaces of arbitrary genus $g$.  However, 
in order to get sharp results  when $g > 0$ 
  we would need detailed information
 about the question of which homology classes are represented by
embedded $J$-holomorphic curves for  {\it arbitrary} (and hence non
generic) tame almost complex structures $J$.  This is a rather delicate matter
that will be explored in a later paper.

The main idea in this paper is to consider the fibration 
$$
\Symp(M, \om)\cap \Diff_0(M) \to \Diff_0(M) \to \Ss_{[\om]}
$$ 
coming from the action of the identity component $ \Diff_0(M)$ of the
group of diffeomorphisms on the space $\Ss_{[\om]}$  of all symplectic
structures on $M$ that are isotopic to $\om$, and then to center the
arguments on the spaces $\Ss_{[\om]}$ rather than on 
 the groups $\Symp(M, \om)$.  This point of view was inspired by
Kronheimer's construction in~\cite{K} of  families of symplectic forms  that
represent elements in $\pi_*(\Ss_{[\om]})$  with nonzero image under the
boundary map  $$
\p: 
\pi_*(\Ss_{[\om]}) \to \pi_{*-1}(\Symp(M, \om)\cap \Diff_0(M)).
$$
This approach allows us to recover some  of the main results in
~\cite{AM}, though  it does not lead to a complete calculation of the
homotopy type of $\Symp(M, \om)$.  So far, the only such calculation besides
Gromov's is by Anjos~\cite{An}, who recently found the homotopy type of the
symplectomorphism group of $S^2\times S^2$ when the ratio $\mu$ of
the sizes of the spheres lies  in the interval $(1, 2]$.

An essential tool in the argument is the existence of many $J$-holomorphic
curves in ruled surfaces, which is a reflection of the fact that they have
many nonzero Seiberg--Witten invariants.  It is not at all clear what can be
said about the topology of $\Symp(M,\om)$ for general $4$-manifolds that
typically have few such curves.  Also, nothing is known in higher dimensions. 
For example, it is not known whether the group of compactly supported
symplectomorphisms of  $\R^{2n}$ is contractible when $n > 2$.

 \subsection{Statement of main results}\label{ss:main}

 Ruled surfaces are compact smooth
$4$-manifolds $M$ that fiber over a Riemann surface $\Si = \Si_g$ of genus
$g$ with fiber $S^2$.  There are two topological types, the product $\Si\times
S^2$ and the total space $M_\Si$ of a nontrivial fibration over $\Si$.  
As in~\cite{AM}, the two cases are analogous.  For simplicity, we will restrict
ourselves here to the case of the trivial fibration.

It is known from the work of 
Taubes, Li--Liu and Lalonde--McDuff (for
detailed references see~\cite{LM2}) that every symplectic form $\om$ on 
$\Si\times S^2$ is  diffeomorphic to a scalar multiple of one of the 
standard forms $\om_\mu$, where
$$
\om_\mu = \mu \si_\Si +  \si_{S^2},\qquad \mu > 0.
$$
 (Here $\si_Y$ denotes an area form on
the Riemann surface $Y$ with total area $1$.)
We denote by $G_\mu = G_\mu^g$ the subgroup  $\Symp(M,
\om_\mu)\cap \Diff_0(M)$ of the group of symplectomorphisms of $(M,
\om_\mu)$.  When $g > 0$ $\mu$ ranges over all positive numbers. 
However, when $g=0$ there is an extra symmetry:  interchanging the two
spheres gives an isomorphism $G_\mu^0 \cong G_{1/\mu}^0$.  It is
important for our arguments that the base sphere be at least as large as the
fiber: see Lemma~\ref{le:inc}.  Hence in this case we take $\mu$ in the range
$[1,\infty)$.

 The groups $G_\mu^g$ were first studied
by Gromov~\cite{G} in the  case when $g=0$.  He showed that 
when $\mu = 1$, i.e. when both spheres have the same size,
$G_1^0$ is connected and
deformation retracts to the Lie group $\SO(3)\times \SO(3)$.  He also
pointed out that as soon as $\mu$ gets bigger than $1$, 
 a new element of infinite order appears in the fundamental
group $\pi_1(G_\mu^0)$.  The key idea in his proof was to look at the action
of $G_\mu^0$ on the contractible space $\Jj_\mu$ of $\om_\mu$-compatible
almost complex structures.  

These ideas were taken much further by Abreu~\cite{A} and
Abreu--McDuff~\cite{AM} and led to a calculation of the rational homotopy
type of both $G_\mu^0$ and $BG_\mu^0$ for all $\mu$.
Three of these results are relevant here.

\begin{prop}\label{prop:am1}  As $\mu\to \infty$,  the groups $G_\mu^g$
tend to a limit $G_\infty^g$ that
has the homotopy type  of the identity component $\Dd_0^g$ of the group 
of fiberwise diffeomorphisms
of $M = \Si_g\times S^2$.
\end{prop}

Here $\Dd_0=\Dd_0^g$ is the identity component of the group $\Dd = \Dd^g$
of all diffeomorphisms $\phi$  that fit into the commutative diagram $$
\begin{array}{ccc} M & \stackrel \phi{\to} &M\\
\downarrow & &\downarrow\\
\Si & \stackrel {\phi'}{\to}  & \Si,
\end{array}
$$
where the vertical arrows are projection to the first factor.  $\Dd_0^g$
can also be
described as the set of all diffeomorphisms that are isotopic to the
identity and preserve the degenerate form $\si_\Si$.  The calculation of the
limit of the $G_\mu^g$ is straightforward: almost the hardest part
is the definition of the maps $G_\mu^g \to G_{\mu + \eps}^g$.  In
\S\ref{ss:bi} below we give an indirect but conceptually simple definition of
these maps.

Observe also that  $\Dd_0^0$ is homotopy equivalent to $\SO(3)\times\SO(3)
\times \Om^2(S^3)$, where $\Om^2(S^3)$ denotes the double loop space of
$S^3$: see \S\ref{ss:dd0}.  
Gromov's new element in $\pi_1(G_\mu^0)$ maps to the generator
of $\pi_1(\Om^2(S^3))$.  Geometrically, it can be represented by the 
``rotation about the diagonal and anti-diagonal", i.e. by the family of maps
$$
S^1\times S^2\times S^2 \to S^2\times S^2:\quad
(t,z,w)\mapsto (z, R_{z,t}(w)),
$$
where $R_{z,t}$ is the rotation of the fiber sphere $S^2$ through the angle
$2\pi t$ about the axis through the point $z\in S^2$.   

Here is the second main result from~\cite{AM}.

\begin{prop}\label{prop:am2}
When $\ell  < \mu \le \ell + 1$ for some integer
$\ell \ge 1$,  
$$
H^*(G_\mu^0, \Q) = \La(t, x, y)\otimes \Q[w_\ell], 
$$
where $\La(t,x,y)$ is an exterior algebra over $\Q$ with generators $t$ of
degree $1$, and $x,y$ of degree $3$ and $\Q[w_\ell]$ is the polynomial
algebra
on a generator $w_\ell$ of degree $4\ell$.
\end{prop}

In the above statement, the generators $x,y$ come from $H^*(\SO(3)\times
\SO(3))$ and $t$ corresponds to Gromov's element in $\pi_1(G_\mu^0), \mu >
1$.   Thus the subalgebra $\La(t, x, y)$ is the pullback of $H^*(\Dd_0^0, \Q)$
under the map $G_\mu \to \Dd_0^0$.    The other generator
$w_\ell$ is fragile, in the sense that it disappears (i.e. becomes null
cohomologous) when $\mu$ increases.

Thus the  map $G_\mu^0 \to \Dd_0^0$
is not a homotopy equivalence for any $\mu$.   
However, as we shall see in Corollary~\ref{cor:d0} the next statement
implies that  it  induces a surjection
on rational homotopy when $\mu> 1$.  

\begin{prop}\label{prop:am3} When $\mu > 1$, $G_\mu^0$ is connected and
$ \pi_1(G_\mu^0) = \pi_1(\Dd_0^0) = \Z + \Z/2\Z + \Z/2\Z.
$ 
\end{prop}

In this paper we give a simplified proof of  Proposition~\ref{prop:am1} above
and a new proof of the following result about the groups
$G_\mu^0$.  Part $(i)$ of the theorem below  sharpens one aspect of 
Proposition~\ref{prop:am2}, while part $(ii)$ is somewhat weaker. 

\begin{thm}\label{thm:main}\begin{itemize}
\item[(i)]
The homotopy type of $G_\mu^0$
is constant on all intervals $(\ell, \ell + 1]$ with $\ell \ge 1$.
Moreover, as $\mu$ passes the integer $\ell + 1 $, $\ell \ge 0$, the groups
$\pi_i(G_\mu^0), i\le 4\ell -1,$ do not change.

\item[(ii)]  There is an 
element  $\rho_{\ell}\in \pi_{4\ell}(G_\mu^0)\otimes \Q$ 
when  $\ell  < \mu \le \ell +1$
 that vanishes for $\mu > \ell +1$. \end{itemize} \end{thm}

Here, the elements $\rho_\ell$ are dual to the classes
$w_\ell\in H^{4\ell}(G_\mu^0, \Q)$ of Proposition~\ref{prop:am2}.  
This result is proved by methods that should in principle work
for base manifolds $\Si$ of any genus.  (In contrast the proofs in~\cite{AM}
use ideas that definitely do not hold when $g > 0$: see
Remark~\ref{rmk:homog}.)  A natural
conjecture would be that the analog of  $(i)$ should
hold for all $g$, perhaps with the proviso that $\ell \ge g$.
Statement $(ii)$ should also hold with the dimension of $\rho_\ell$ being
taken as $4\ell + 2g$.   
As we shall see in \S\ref{ss:pr} below, what we
would need to  extend our results  in this way is a better knowledge of
which classes in $H_2(M)$ can be represented by embedded
$J$-holomorphic curves.  We will content ourselves here with stating a few
easy partial results.
The first result concerns the case $g = 1$ and the second
is an analog of Proposition~\ref{prop:am3}.

\begin{prop}\label{prop:g1}\begin{itemize}
\item[(i)]
The homotopy type of the groups $G_\mu^1$ is constant for $0 < \mu \le 1$.
Moreover,  the  map $G_\mu^1 \to \Dd_0^1$ gives rise to 
isomorphisms on $\pi_i, i = 0,1,$  for all $\mu > 0$, and  isomorphisms
on $\pi_i, i = 2, \dots, 5,$ for all $\mu > 3/2$.

\item[(ii)]
There is a nonzero element $\rho\in \pi_2(G_1^1)$ that vanishes in
$\pi_2(G_\mu^1)$ when $\mu > 1$.
\end{itemize}
\end{prop}

\begin{prop}\label{prop:infty} \begin{itemize}
\item[(i)]   If $g = 2k$ or $2k+1$ and $\mu > k$ then
the map $G_\mu^g \to \Dd_0^g$ gives rise to an 
isomorphism on $\pi_i$ for $i \le 2g -1$. In particular, $G_\mu^g$ is
connected whenever $\mu > k$.
\item[(ii)]  When $g> 0$, the induced map
$G_\mu^g \to \Dd_0^g$ induces a surjection on rational homotopy groups for
all $\mu > 0$.  
\end{itemize}
\end{prop}

\NI
{\bf Further Questions}\MS

\NI$\bullet$  As noted in Proposition~\ref{prop:g1} above, 
the homotopy type of $G_\mu^1$ is constant for $0 < \mu \le 1$.  It would be
interesting to figure out its structure.  In the genus $0$ case, it was
shown in~\cite{AM} that there is a family of 
Lie subgroups $K_k$ of $\Diff_0(M)$ such that
the whole rational homotopy type of $G_\mu^0$ is generated by the
 $K_k$ for $k < \mu$.\footnote
{
By Anjos~\cite{An}, this statement remains true for the integral homotopy
type of $G_\mu^0, 1 < \mu \le 2$.}  
In particular, the new fragile elements $\rho_k$ are higher
Whitehead (or Samelson) products of elements in $\pi_*(K_k)$.  It is not
clear if there is an  analogous  result when $g= 1$.  In particular, it would be
nice to find a specific representative of the element $\rho$ in
$\pi_2(G_1^1)$. \MS

\NI$\bullet$  It seems very likely that the group $G_\mu^g$ is always
connected.  Equivalently, it is likely that the map
$$
\pi_0(\Symp(M,\om_\mu)) \to \pi_0(\Dd^g),
$$
which can be shown to be a surjection as in the proof of 
Proposition~\ref{prop:infty} $(ii)$,
is in fact an isomorphism.
\MS

\NI$\bullet$ 
Another direction in which one might speculate concerns the relation of the
groups $G_\mu^g$  to the full group of diffeomorphisms.  For example,
consider the case  $M = S^2\times S^2$.  Then $\Diff_0(M)$ contains two
copies of $\Dd_0^0$ which we might call $\Dd_{{\it left}}$ and $\Dd_{{\it
right}}$ corresponding to the two projections $M\to S^2$.  They intersect in
$\SO(3)\times \SO(3) = G_1^0$.   It would be interesting to know how much of
the topology of $\Diff_0(S^2\times S^2)$ is generated by these subgroups.
Could is possibly be true that, just as the groups $K_k$ generate the
rational homotopy of $G_\mu^0$, these two groups generate the rational
homotopy of $\Diff_0(S^2\times S^2)$?
 In particular, what is the Samelson product of the two elements
$\tau_{{\it left}}, \tau_{{\it right}}$ of $\pi_1(\Diff_0(S^2\times S^2))$ given by the
generator $\tau$ of $\pi_1(\Dd_0^0)$?
There are similar questions when $g > 0$.
However there  is less
topology to play with, since the groups $\Diff_0(\Si_g)$ are  contractible
when $g > 1$. \MS

{\bf  Organization of the paper}
In \S\ref{s:out} we first show how to define and calculate
the limit $G_\infty$, and  then  describe the main
points in the proofs of Theorem~\ref{thm:main},
Proposition~\ref{prop:g1} and Proposition~\ref{prop:infty}  $(i)$.  The
next two sections discuss some  needed techniques:  symplectic inflation 
in~\S\ref{s:infl} and the construction of embedded curves in~\S\ref{s:emb}.  
Finally, in \S\ref{s:dd} we study the topology of the fiberwise groups
$\Dd_g^0$ and prove  Proposition~\ref{prop:infty} $(ii)$.  This
 is elementary except for the statement that
$\pi_1(G_\mu^g)$ maps onto $\pi_1(\Dd_0^g)$ when $\mu \le 1$ and $ g >
1.$

\section{Outline of proofs}\label{s:out}

In this section we define and calculate
the limit $G_\infty$ and then show how to understand the
relationship between the different spaces $\Aa_\mu$. 

\subsection{The basic idea}\label{ss:bi}

 There is no very direct map $G_\mu \to
G_{\mu + \eps}$ for $\eps > 0$, and the proof 
of Proposition~\ref{prop:am1} given in~\cite{AM} is rather
clumsy.    This section gives a
streamlined version of the argument that is based on looking at the
fibration 
$$ 
G_\mu \to \Diff_0(M) \to \Ss_\mu.
$$
Here $\Diff_0(M)$ is the identity component of the group of
diffeomorphisms of $M$ and $\Ss_\mu$ is the space of all symplectic forms
on $M$ in the cohomology class $[\om_\mu]$ that are isotopic to $\om_\mu$.
By Moser's theorem the group $\Diff_0(M)$ acts transitively on $\Ss_\mu$ 
via an action that we will write as 
$$
\phi\cdot \om = \phi_*(\om) = (\phi^{-1})^* \om,
$$
so that
$\Ss_\mu$ is simply the homogeneous  space $\Diff_0(M) / G_\mu$.  
Let $\Aa_\mu$ denote the space   of
almost complex structures that are tamed by some form in
$\Ss_\mu$.\footnote {
Recall that $\om$ is said to tame $J$ if $\om(v, Jv) > 0$ for all $v\ne 0$.}  

\begin{lemma}\label{le:sa} $\Ss_\mu$ is homotopy equivalent to $\Aa_\mu$.
\end{lemma}
\proof{} Let $\Xx_\mu$ be the space of pairs
$$
\Xx_\mu= \{(\om,J)\in \Ss_\mu\times \Aa_\mu:
\om\mbox{ tames } J\}. 
$$ 
Then both projection maps $\Xx_\mu\to 
\Aa_\mu, \Xx_\mu\to \Ss_\mu$ are fibrations with contractible fibers
(see~\cite{MS2} Ch 2.5), and
so are homotopy equivalences.\QED

 \begin{lemma}\label{le:inc}  If $g  > 0$ then
$\Aa_\mu\subset \Aa_{\mu + \eps}$ for all $\mu, \eps > 0$. When $g=0$ this
holds provided $\mu \ge 1$. \end{lemma} 
\NI{\bf Sketch of proof:\,\,}  It is well known that for each $J\in
\Aa_\mu$, $M$ admits a 
 foliation whose leaves are $J$-holomorphic spheres  representing the
fiber class $F = [pt \times S^2]$.  For the sake of completeness, we sketch a
proof in Lemma~\ref{le:fib} below.  Observe that this is where we use the
fact that $\mu \ge 1$ when $\Si = S^2$:  if the base is smaller than the fiber,
the fiber class can degenerate so that  $\pi_J$ may not exist.

It follows that there is a projection $\pi_J:M\to \Si$ 
onto the leaf space $\Si$ of this foliation.  We show in~\S\ref{s:infl}  that if 
$\si$ is any area form on $\Si$ then its pullback $\pi_J^*(\si)$ is
$J$-semi-tame, i.e 
$$
\pi_J^*(\si)(v, Jv) \ge 0,\quad v\in TM.
$$
Granted this, if $\om\in \Ss_\mu$ tames $J$, then $\om + \ka \pi_J^*(\si)$
also tames $J$ for all $\ka > 0$.   The result is then immediate.\QED

\begin{cor} For any $\mu, \eps, \eps' > 0$ there are  maps
$\Ss_\mu \to \Ss_{\mu + \eps}$ and $G_\mu \to G_{\mu + \eps}$ that 
are well defined up to homotopy and make
the following diagrams homotopy commute: 
$$
\begin{array}{ccccccc}
{(a)} & & G_{\mu} &\to & \Diff_0(M) & \to & \Ss_{\mu}\\
& & \downarrow & & \;\downarrow = & & \downarrow\\
& & G_{\mu + \eps} &\to & \Diff_0(M) & \to & \Ss_{\mu + \eps},\\
& & & & & & \\
{ (b)} & & 
G_{\mu} &\to & G_{\mu+\eps} & & \\
& & & \searrow &\downarrow & & \\
& & && G_{\mu + \eps+\eps'}& &.
\end{array}
$$
\end{cor}
\proof{}  The maps $\Ss_\mu \to \Ss_{\mu + \eps}$ are induced from
the inclusions $\Aa_\mu \subset \Aa_{\mu + \eps}$ using the homotopy
equivalences $\Ss_\mu \simeq \Aa_\mu$ in Lemma~\ref{le:sa} above.  
Since $G_\mu$ is the fiber of the map
$\Diff_0(M) \to \Ss_{\mu}$, there are induced maps $G_\mu \to 
G_{\mu + \eps}$ making diagram $ (a)$ homotopy commute.  The rest is
obvious.\QED

This  corollary illustrates the essential feature
of our approach.  Statements that are true only up to homotopy on the level
of the groups $G_\mu$ are true on the nose for the spaces $\Aa_\mu$.
The discussion below shows that it is easy to understand what happens to
the groups $G_\mu$ as $\mu$ increases.  Much of the rest of the paper is
devoted to understanding (on the level of the spaces $\Aa_\mu$) what
happens as $\mu$ decreases.  For this we use the Lalonde--McDuff
technique of symplectic inflation. 
\MS

\NI
{\bf Proof of Proposition~\ref{prop:am1}.}

We first show that we can understand the limit $G_\infty = \lim G_\mu$
by studying the space $\cup_\mu \Aa_\mu$.  
Let $J_{split}$ be the standard product almost complex structure on $M$.
Because $J_{split}$ is tamed by $\om_\mu$ the map $\Diff_0(M) \to
\Ss_\mu$ lifts to
$$
\Diff_0(M) \to\Xx_\mu:\quad \phi\mapsto (\phi_*(\om_\mu),
\phi_*(J_{split})).
$$
Composing with the projection to $\Aa_\mu$ we get a map
$$
\Diff_0(M) \to\Aa_\mu:\quad \phi\mapsto 
\phi_*(J_{split})
$$
that is not a fibration but has homotopy fiber $G_\mu$.
Since $\Aa_\mu$ is an open subset of $\Aa_{\mu + \eps}$ for all $\eps > 0$,
the homotopy limit $\lim_\mu \Aa_\mu$ of the spaces $\Aa_\mu$ is
homotopy equivalent to the union $\Aa_\infty = \cup_\mu \Aa_\mu$.  
Hence $G_\infty$, which is defined to be the homotopy limit of the $G_\mu$,
is  homotopy equivalent to the homotopy fiber of the map
$\Diff_0(M) \to \Aa_\infty$.

To understand $\Aa_\infty$ we proceed as follows. 
Let $\Fol$ be the space of all smooth foliations of $\Si\times S^2$ by
spheres in the fiber class $F = [pt \times S^2]$.  Since $S^2$ is compact and
 simply connected, each leaf of this foliation has trivial holonomy and
hence has a neighborhood that is diffeomorphic to the product $D^2\times
S^2$ equipped with the trivial foliation with leaves $pt\times S^2$.  It follows
that $\Diff(M)$ acts transitively on 
$\Fol$ via the map $\phi\mapsto \phi(\Ff_{split})$, where $\Ff_{split}$ is the
flat foliation  by the spheres $pt \times S^2$.  Similarly, $\Diff_0(M)$
acts transitively on the connected component $\Fol_0$ of $\Fol$ that
contains $\Ff_{split}$.  Hence  there is a fibration sequence 
$$ 
\Dd\cap \Diff_0(M) \to
\Diff_0(M) \to \Fol_0. 
$$
It is not hard to see that the group $\Dd\cap \Diff_0(M)$ is connected,
and so equal to $\Dd_0$: see Corollary~\ref{cor:con}.

Next, observe that there is a map $\Aa_\infty \to \Fol_0$ given by taking
$J$ to the foliation of $M$ by $J$-spheres in class $F$.  Standard arguments
(see for example~\cite{MS2} Ch 2.5) show that this map is  a
fibration with contractible fibers.  Hence it is a homotopy equivalence.
Moreover, it fits into the
commutative diagram:
$$
\begin{array}{ccc}
  \Diff_0(M)  & \to  & \Aa_\infty\\
\downarrow & & \downarrow\\
\Diff_0(M) & \to & \Fol_0,
\end{array}
$$
where the map $ \Diff_0(M)  \to   \Aa_\infty$ is given as above by the action
$\phi\mapsto \phi_*(J_{split})$.  Hence there is an induced homotopy
equivalence from the homotopy fiber $G_\infty$ of the top row to the fiber
$\Dd_0$ of the second.
  \QED

\begin{rmk}\rm  Implicit in the above argument is the following
 description of the map $G_\infty \to \Dd_0$.  Let $\Jj_\mu$ denote the
 space of all almost complex structures tamed by $\om_\mu$.  
Since the image of the group $G_\mu$ under the map $\Diff_0(M) \to
\Aa_\mu$ is contained in $\Jj_\mu$ there is
a commutative diagram
$$
\begin{array}{ccc} & & \Dd_0\\
&  & \downarrow\\
G_\mu &\stackrel{\io}\longrightarrow & \Diff_0(M)\\
\downarrow & & \downarrow\\
\Jj_\mu & \longrightarrow & \Fol.
\end{array}
$$
Because $\Jj_\mu$ is contractible, the inclusion $\io:
G_\mu \to \Diff_0(M)$  lifts to a map  $
\Tilde\io: G_\mu \to \Dd_0.$  Now take the limit to get $G_\infty \to \Dd_0$.
\end{rmk}

\subsection{The stratification of $\Aa_\mu$}\label{ss:pr}

The following definition is the key to understanding the relation between
the  $\Aa_\mu$ as $\mu $ decreases.

\begin{defn}\label{def:k0}\rm
Given $\ell > 0$, let $\Aa_{\mu,\ell}$ be the subset of $J\in \Aa_\mu$
consisting of elements that admit a  $J$-holomorphic curve in class
$A - \ell F$.   Here $F$ is the fiber class as before, and $A = [\Si \times pt]$.
Further, define
$$
\Aa_{\mu, 0} = \Aa_\mu - \cup_{\ell > 0} \Aa_{\mu,\ell}.
$$
\end{defn}

Since the fiber class $F$ is always represented and since distinct
$J$-holomorphic curves always intersect positively,\footnote
{See~\cite{LM1} for an
outline of the main points of the theory of $J$-holomorphic curves,
 and Sikorav~\cite{Si} for
detailed proofs of positivity of intersections.}
 the class $pA + qF$ has a $J$-holomorphic representative only if 
$$
p = (pA + qF)\cdot F \ge 0.
$$
Thus $\Aa_{\mu, 0}$ can be characterized as the set of $J$ for which
the only classes with $J$-holomorphic representatives are 
$pA + qF$ with $p,q\ge 0$. A similar argument shows that the sets
$\Aa_{\mu, \ell}$ are disjoint for $\ell\ge 0$.
Note also that $\Aa_{\mu, \ell}$ is nonempty only if $\ell < \mu$:  if $J$
admits  a $J$-holomorphic curve in class
$A - \ell F$, any form $\om\in \Ss_\mu$ that tames $J$ must be positive on
this curve.  Hence $\om(A - \ell F) = \mu - \ell > 0$.  

The next lemma
shows that the subsets $\Aa_{\mu, \ell}$ are well behaved. 

\begin{lemma}\label{le:str} For all $\ell > 0$, $\Aa_{\mu, \ell}$ is a
Fr\'ech\'et suborbifold of $\Aa_\mu$ of codimension $n(\ell)= 4\ell- 2 +
2g$. 
\end{lemma}
\proof{}
 This is proved in the genus $0$ case by Abreu~\cite{A}, and the proof
for general $g$ is similar.  (See also~\cite{M3}.)  The idea is  the following. 
Let $\Mm$ denote the moduli space of all triples $(u,j,J)$
where $J \in \Aa_\mu$, $(\Si, j)$ is a marked Riemann surface 
(so $j $ is in Teichm\"uller space) and  
 $u:(\Si, j)\to M$ is a $J$-holomorphic curve in class 
$A - \ell F$.   By standard theory $\Mm$ is a Fr\'ech\'et manifold.
The forgetful map
$$
\Mm \to \Aa_\mu,\quad  (u,J)\mapsto J,
$$
has image $\Aa_{\mu, \ell}$.  Moreover its fiber 
at $J$ is  the fiber of the projection from Teichm\"uller space to the moduli
space of Riemann surfaces of genus $g$.  Hence it varies by at most a finite
group.   This gives  $\Aa_{\mu, \ell}$ the structure of a Fr\'ech\'et orbifold.

It remains to check that it has finite codimension in $\Aa_\mu$.
To see this, note that the curve $C= {\rm Im}\, u$
 is embedded since each fiber intersects it exactly
once, positively.  (By positivity of intersections, if $u$ were  singular at
$z \in \Si$, then the intersection number of $C$ with the fiber
through $u(z)$ would have to be $> 1$.)  Moreover 
$c_1(A - \ell F) < 0$.  Hofer--Lizan--Sikorav
showed in~\cite{HLS} that in this case the linearization 
$$
Du: \;\;C^\infty(u^*TM) \to  \Om_J^{0,1}(u^*TM)
$$
of the
generalized Cauchy--Riemann operator $\p_J$ at $u$ has cokernel of
constant  rank equal to $n(\ell)$.  Moreover,
Abreu's arguments show that this cokernel can be identified with the normal
space to $\Aa_{\mu, \ell}$ at $J$ as follows.  Since
$\Aa_\mu$ is an open subset in the space of all almost complex structures on
$M$, the tangent space $T_J\Aa_\mu$ consists of all smooth sections $Y$ of
the bundle $\End_J(TM)$ of endomorphisms of $TM$ that anticommute with
$J$.  The claim is that the restriction map  
$$
\End_J(TM) \to \Om_J^{0,1}(u^*(TM)):\quad Y\mapsto Y\circ du\circ j
$$
identifies the fiber of the normal bundle at $J$ with the cokernel of $Du$.
By standard arguments (see for example Proposition~3.4.1 in~\cite{MS1})
we can restrict attention here to sections $Y$ with support in some open set
that intersects the curve $C$. \QED

 The following key result is proved in \S\ref{s:infl} using
the technique of symplectic inflation from~\cite{LM2}.   

\begin{prop}\label{prop:key}   Let $J \in \Aa_\mu$.  If the class $pA + qF$ is
represented by an embedded  $J$-holomorphic curve $C$ for some $p > 0,
q\ge 0$, then  $J\in A_{\eps + q/p}$ for all $\eps > 0$.
\end{prop}

The idea of the proof is to move the cohomology class of the taming form
$\om$ for $J$ towards the Poincar\'e dual of $[C]$ by adding to $\om$ a large
multiple of a form $\tau_C$ that represents this class: see
\S~\ref{s:infl}. It follows that if
we can find classes of small ``slope" $q/p$ that have  embedded
$J$-holomorphic representatives, we can decrease the  parameter $\mu$ of
the set $\Aa_{\mu}$ containing $J$. Note that if  $J\in \Aa_{\mu, \ell}$ and
$pA + qF$ has an embedded  $J$-representative, then positivity of
intersections implies that $$
(A - \ell F)\cdot (pA + qF) \ge 0, \quad \mbox{i.e.} \;\;\frac qp \ge \ell.
$$
Thus if  $J\in \Aa_{\mu, \ell}$
we cannot hope to  represent a class  of smaller slope than $A + \ell F$.

The following result is proved in~\S\ref{s:emb}.

\begin{lemma}\label{le:g0}
If the base
has  genus $g= 0$ and $J\in \Aa_{\mu, \ell}$ then
the class $A + \ell F$  always has an embedded 
$J$-representative. 
\end{lemma}

\begin{cor}\label{cor:g0}  When $g = 0$, $A_{\mu, \ell} = \Aa_{\ell+\eps,
\ell}$ for all $\eps > 0$.  Thus the sets 
$
\Aa_\mu = \cup_{\ell < \mu} \Aa_{\mu, \ell}$ are constant on all intervals
$k < \mu \le k+1$.
\end{cor}

\NI
{\bf Proof of Theorem~\ref{thm:main} $(i)$}

We must show that the homotopy type of the groups $G_\mu^0$ is constant
for $\mu \in (\ell, \ell +1]$.  
 Because of the fibration 
$$
 G_\mu\to \Diff_0(M) \to \Ss_\mu,
$$
it suffices to prove this for $\Ss_\mu$ and hence, by Lemma~\ref{le:sa}, for
$\Aa_\mu.$   But this holds by the preceding corollary.

To prove the second statement, use Lemma~\ref{le:str} which says that when
$\mu$ passes $\ell + 1$ the topology of $\Aa_\mu$ changes by the addition
of a stratum of codimension $n(\ell) = 4\ell + 2.$   Hence its homotopy groups
$\pi_i$ for $i \le 4\ell$ do not change.  Thus the groups $\pi_i(G_\mu^0)$ for
$i\le 4\ell -1$, do not change. \QED

It is not clear whether Lemma~\ref{le:g0} continues to hold for arbitrary $g$. 
The problem is that we  have to understand the $J$-holomorphic curves for
{\em every} $J$ not just the generic ones, or even just the integrable ones.
In \S\ref{s:emb} we exhibit the existence of enough embedded
$J$-holomorphic curves to prove:

\begin{lemma}\label{le:g}  If $g = 2k$ or $2k+1$ for some $k\ge 0$,
the space $A_{\mu, 0}^g$ is constant for  $\mu > k$.  Moreover,
 the space $A_{\mu, 1}^1$ is constant for $\mu > 3/2$.
\end{lemma}

\begin{cor}  Part (i) of Propositions~\ref{prop:g1} and~\ref{prop:infty} hold.
\end{cor}
\proof{}  Because  $\Aa_\mu^g =
\Aa_{\mu, 0}^g$ when $\mu \le 1$, 
Lemma~\ref{le:g} implies that  $\Aa_\mu^1$ is constant for
$\mu\in (0,1]$.   Hence the first statement of
Proposition~\ref{prop:g1} $(i)$ holds by the arguments given  above. 
Moreover, if $\mu > 3/2$ and $\eps > 0$
the difference $\Aa_{\mu + \eps}^1 - \Aa_\mu^1$ is contained in a union of
strata of codimension at least $n(2) = 8$.  Therefore, the inclusion
$G_\mu^1\to G_{\mu + \eps}^1$ induces an isomorphism on $\pi_k$ for $k \le
5$ and hence an isomorphism
$$
\pi_k(G_\mu^1) \stackrel{\cong}\to \pi_k(\lim G_\mu^1) = \pi_k(\Dd_0^1)
$$
for these $k$.  The proof of
 Proposition~\ref{prop:infty} $(i)$ is similar. \QED

Part $(ii)$ of these propositions are proved by different methods.  
Proposition~\ref{prop:g1} $(ii)$ is proved together with
Theorem~\ref{thm:main} $(ii)$  in the next subsection.   By
contrast Proposition~\ref{prop:infty} $(ii)$ is proved essentially by direct
calculation: see \S\ref{ss:dd}.

\subsection{Topological changes as $\mu$ passes an integer}\label{ss:frag}

 We will write $\Xx_{\mu,k}$ for the inverse
image of $\Aa_{\mu,k}$ in $\Xx_\mu$.  Thus $\Xx_{\mu, k}$ consists of
all pairs $(\om, J) \in \Xx_\mu$ for which $J \in \Aa_{\mu, k}$.
 Note that
whenever  $\Xx_{\mu, k}$ is nonempty the map $\Xx_{\mu, k} \to
\Ss_\mu$ is a fibration whose fiber will be denoted $\Jj_{\mu, k}$. 
Thus, the contractible space $\Jj_\mu$ of almost complex structures
tamed by  $\om_\mu$ is the union of the $\Jj_{\mu, k}$.  Since 
$\Jj_\mu$ is open in $\Aa(\om_\mu)$,
 each $\Jj_{\mu, k} = \Jj_\mu\cap \Aa_{\mu, k}$ is a
suborbifold of  $\Jj_\mu$, 

As $\mu$ passes the integer
$k$, a new stratum of codimension $n = n(k) = 4k - 2 + 2g$ appears in the
spaces $\Aa_\mu, \Xx_\mu, \Jj_\mu$.  Since $\Jj_\mu$ is contractible
and $\Jj_{\mu, k}$ is an orbifold, there is an isomorphism
of rational homology
$$
H_i(\Jj_{\mu, k}) \cong H_{i+n-1}(\Jj_{\mu}-\Jj_{\mu, k}).
$$
In particular, there is a nonzero element $\al\in H_{n-1}(\Jj_{\mu}-\Jj_{\mu,
k})$ that is generated by the sphere linking the new
stratum.  Since $\al$ is spherical it must represent a nonzero class in
the  homotopy group $\pi_{n-1}(\Jj_{\mu}-\Jj_{\mu, k})$.
When $g = 0$ these elements turn out to generate the whole of the
topology change in $G_\mu$ as $\mu$ passes $k$.   It is not
clear whether this
extends  to higher genus.  However, as we now see, in some circumstances it
is possible to work out the effect of these new elements $\al$ on the rational
homotopy of $G_\mu$.

\begin{prop}  Let $n = 4k - 2 + 2g$ and suppose that
 the inclusion $\Aa_k^g \to \Aa_{k + 1}^g - \Aa_{k + 1,k}^g$
induces an isomorphism on $\pi_i$ for all $i \le n$.  Suppose further 
that the map $G_\mu^g \to \Dd_0^g$ induces an isomorphism on 
$\pi_{n-1}\otimes \Q$ for $\mu = k+1$ and is surjective for $\mu = k$.  Then
there is a nonzero element $\rho \in \pi_{n-2}(G_k)\otimes \Q$  that
vanishes in $\pi_{n-2}(G_\mu)\otimes \Q$ when $\mu > k$.
\end{prop}
\proof{}
Consider the following diagram: $$
\begin{array}{ccccc}
\Jj_{k+1}-\Jj_{k+1, k} &\to & 
\Xx_{k+1}-\Xx_{k+1, k} &\to & \Ss_{k+1}\\
\downarrow & & \downarrow & & \downarrow \\
\Jj_{k+1} &\to & \Xx_{k+1} &\to & \Ss_{k+1}.
\end{array}
$$
The rows here are fibrations.  As already noted, 
$$
 H_{n-1}(\Jj_{k+1}-\Jj_{k+1, k})\cong  H_0(\Jj_{k+1, k}) \cong \Q,
$$
and is generated by a spherical element
$\al\in \pi_{n-1}(\Jj_{k+1}-\Jj_{k+1, k})$.
Suppose first that $\al$ goes to
$0$ in  $\Xx_{k+1}-\Xx_{k+1, k} $ and so comes from a nonzero element
$\be\in \pi_n(\Ss_{k+1})\otimes \Q$, and consider the fibration
$$
 G_{k+1}\to \Diff_0(M) \to \Ss_{k+1}.
$$
The first
claim is that $\be$ cannot lift to an element $\be'\in
\pi_n(\Diff_0(M))\otimes \Q$.  For if it did,  because $\Diff_0(M)$ also acts on
$\Ss_k$, $\be$ would be homotopic to an element in  
$$
\pi_n(\Ss_k)\otimes \Q = \pi_n(\Xx_k)\otimes \Q
 = \pi_n(\Xx_{k+1}-\Xx_{k+1, k})\otimes \Q
$$
and so would go to zero under the boundary map from
$\pi_n(\Ss_{k+1})$ into $\pi_{n-1} (\Jj_{k+1}-\Jj_{k+1, k})\otimes \Q$.  
This contradicts the fact that it has nonzero
image $\al$ under this map.

Hence $\be$ must give rise to a nonzero element in
$\pi_{n-1}(G_{k+1})\otimes \Q$ that goes to zero in $\Diff_0(M)$. 
But there are no such elements:
 by  hypothesis, the map $\pi_{n-1}(G_{k+1})\otimes \Q \to 
\pi_{n-1}(\Dd_0^g)\otimes \Q$ is an isomorphism, and 
it follows from the explicit calculation of the groups
$\pi_{i}(\Dd_0^g)\otimes \Q$ in \S\ref{ss:dd0} that 
$\pi_*(\Dd_0^g)\otimes \Q$ injects into $\pi_*(\Diff_0(M))$.
 Hence the original hypothesis must be wrong, i.e. $\al$ must map to a
nonzero element $\al'$ in $\pi_{n-1}(\Xx_{{k+1}}-\Xx_{{k+1}, k}) \otimes \Q = 
\pi_{n-1}(\Xx_{k}) \otimes \Q$

There are now two possibilities:\MS

\NI
$(i)$  The element $\al'\in \pi_{n-1}(\Xx_k)\otimes \Q$ comes from
$\ga\in \pi_{n-1}(\Diff_0(M))\otimes \Q$; or\MS

\NI
$(ii)$  $\al'$ gives rise to a nonzero element $\al''\in \pi_{n-2}(G_k)
\otimes \Q$.\MS

Note that these elements $\al', \al''$ are fragile in the sense that they vanish
in $\Xx_\mu$ and $G_{\mu}$ when $\mu > k$  since $\al$ does.  Thus in case
$(i)$ the element  $\ga\in \pi_{n-1}(\Diff_0(M))\otimes \Q$ is in the image
of  $\pi_{n-1}(G_\mu)\otimes \Q$ for $\mu > k$ but 
is not in this image for $\mu =
k$.  Since this contradicts our hypothesis, we must be in case $(ii)$. 
This completes the proof.
\QED

\NI
{\bf Proof of Theorem~\ref{thm:main} $(ii)$ and Proposition~\ref{prop:g1}
$(ii)$}

Apply the above proposition with $k = \ell +1$ and $g=0$, so that $n = 4k-2$.
The first hypothesis above is satisfied for all $k$ since
$\Aa_k^0 = \Aa_{k + 1}^0 - \Aa_{k + 1,k}^0$  by Corollary~\ref{cor:g0}.
The second hypothesis holds when $k\ge 2$ by
Theorem~\ref{thm:main} $(i)$: indeed the lowest homotopy group 
of $G_\mu^0$ that
changes   when $\mu> k$ increases has dimension 
$4k > n -1$.  The third
hypothesis holds for $k \ge 2$ by Corollary~\ref{cor:d0}.  This proves
 Theorem~\ref{thm:main} $(ii)$.  

  The proof of Proposition~\ref{prop:g1} $(ii)$ is similar.  
Take $k = 1$.  Then $\Aa_1^1 = \Aa_{2,0}^1$ so that the
first hypothesis holds.  The map $\pi_3(G_\mu^1)\otimes \Q\to
\pi_3(\Dd_0^1)\otimes \Q$ is surjective for all $\mu$  by
Proposition~\ref{prop:infty} $(ii)$, and is an isomorphism for $\mu = 2$ by
Proposition~\ref{prop:g1} $(i)$. Note that the results in \S\ref{s:dd} 
are independent of the current arguments and so it is permissible to use
them here.\QED

\begin{rmk}\rm
It was proved in~\cite{AM} that the elements $\rho_\ell\in
\pi_{4\ell}(G_\mu^0)$ exist only when $\mu\in (\ell, \ell + 1]$.
To prove this in the present context one needs to see that  
the link $\al\in H_*(\Jj_\mu - \Jj_{\mu, \ell})$ cannot be represented in
$\Jj_\mu - \Jj_{\mu, \ell} - \Jj_{\mu, \ell- 1}$, i.e. that its intersection with
the stratum $\Jj_{\mu, \ell-1}$ is homologically nontrivial.  Thus we need
information on the homology of the stratum $\Jj_{\mu, \ell-1}$.
When the base is a sphere we have this information because the strata are
homogeneous spaces of $G_\mu^0$.  The intersection of $\al$ with 
$\Jj_{\mu, \ell-1}$ is just the link of $\Jj_{\mu, \ell}$ in $\Jj_{\mu, \ell-1}$,
and the fact that this is homologically essential is a key step in the
argument:  see~\cite{AM} Corollary~3.2 and its use in the proof of
Proposition~3.7.
 \end{rmk}

\section{Inflation}\label{s:infl}

This section is devoted to the proofs of Lemma~\ref{le:inc}
and Proposition~\ref{prop:key}. \MS

\NI
{\bf Proof of Lemma~\ref{le:inc}.}\,\,  It remains to prove the claim, i.e. that
 if  $\si$ is any
area form on $\Si$ then its pullback $\pi_J^*(\si)$ is $J$-semi-tame.  To see
this,
fix a point $p\in M$ and choose (positively oriented) coordinates $(x_1, x_2,
y_1, y_2)$ near $p$ so that the fibers near $p$ have equation $x_i = const$,
for $i = 1,2,$ and so that the symplectic orthogonal to the fiber $F_p$ at $p$
$$ 
Hor_p =  \{u: \om(u,v) = 0, \mbox{ for all }\; v\in T_p(F_p)\}
$$
is tangent at $p$ to the surface $y_1=y_2 =0$.  Then, at $p$, the form 
$\om_p$ may be written
$$
\om_p = a dx_1\wedge dx_2 + b dy_1 \wedge dy_2
$$
for some constants $a,b > 0$.  Moreover, since $J_p$ preserves the fibers, we
can assume that $J_p$ acts on $T_pM$ via the lower triangular matrix
$$
J_p = \left(\begin{array}{cc} A & 0\\C & J_0\end{array}\right),\quad
\mbox{where}\quad  J_0 = \left(\begin{array}{cc} 0 & -1\\1 & 0
\end{array}\right).
$$
If $u\in Hor_p$, then 
$$
\om_p(u, Ju) = a dx_1 \wedge dx_2 (u, Au) > 0,\quad \mbox{if }\;\; u\ne 0,
$$
because $\om_p$ tames $J_p$.
But $\pi_J^*(\si)$ is just a positive multiple of 
$dx_1 \wedge dx_2$ at $p$.  The claim follows.\QED

To prove Proposition~\ref{prop:key}, it clearly suffices to prove the following
result.  (A weaker version, in which the symplectic forms $\tau_\la$ tame a
suitable perturbation of $J$, is proved in~\cite{LM2}.)

\begin{lemma}[Inflation]\label{le:infl}  Let $J$ be an $\tau$-tame almost
complex structure on a symplectic $4$-manifold $(M, \tau_0)$ that
admits a $J$-holomorphic curve $Z$ with $Z\cdot Z \ge 0$. Then there is a
family $\tau_\la, \la \ge 0,$ of symplectic forms  that all tame $J$ and have
cohomology class $[\tau_\la] = [\tau_0] + \la a_Z$, where $a_Z$ is Poincar\'e
dual to $[Z]$.
\end{lemma}

 We prove this first in the case
when  $Z$ has trivial normal bundle.

\begin{lemma}
Suppose $\tau_0 = \si_Z + dy_1\wedge dy_2$ in $Z\times D^2$
where $(y_1,y_2)$ are rectilinear coordinates on the unit disc $D^2$, and let
$J$ be any $\tau_0$-tame almost complex structure on $Z\times D^2$ such
that $Z\times \{0\}$ is $J$-holomorphic.  Then there is a smooth family of
symplectic forms  $\tau_\la, \la \ge 0,$ on $Z\times D^2$ such that
\begin{itemize}
\item[(i)]
each $\tau_\la$ tames $J$;
\item[(ii)]  for all $\la > 0$,
 $$\int_{pt \times D} \tau_\la = \pi + \la;$$
\item[(iii)]  $\tau_\la = \tau_0$ near $Z\times \p D^2$.
\end{itemize}
\end{lemma}
\proof{}  We will construct $\tau_\la$ to have the form 
$$
\tau_\la = \si_Z + f_\la(r) dy_1\wedge dy_2,
$$
where $f_\la$ is a function of the polar radius $r = \sqrt{y_1^2 +
y_2^2}$ that is chosen to be $\ge 1$ everywhere and $= 1$ near $r = 1$.
The only problem is to make sure that $\tau_\la$  tames $J$. 

 Consider a
point $p = (z, y_1, y_2) \in Z\times D^2$.  We will write elements of
$T_p(Z\times D^2)$ as pairs $(u, v)$ where $u$ is tangent to $Z$ and $v$ is
tangent to $D^2$, and will choose these coordinates so that 
$$
\tau_\la((u,v), (u',v')) =  u^T J_0^T u' + f_\la v^T J_0^T v',
$$
where $^T$ denotes transpose and  $J_0$ is as before.
With respect to the obvious product coordinates, we can
write $J_p$ in block diagonal form $$
J_p = \left(\begin{array}{cc} A & B\\ C & D\end{array}\right)
$$
where $A, B, C, D$ are $2\times 2$ matrices.  
Therefore, we have to choose $f = f_\la(r)$ so  that
$$
\tau_\la((u,v), J_p(u,v))  = u^T J_0^T A u +  u^T J_0^T Bv + f v^T J_0^T C u
+ f v^T J_0^T D v \ge 0 $$
for all $(u,v)\in T_p(Z\times D^2$.
 Because $B = C = 0$ when $r = 0$ and $\tau_0$ tames $J$, there are 
constants $c, c'$ so that the following estimates  hold in some neighborhood
$\{ r\le r_0\}$ of $Z\times \{0\}$: 
\begin{eqnarray*}
\|u\|^2 \le c\,u^T J_0^T A u, & &  \|v\|^2 \le c\,v^T J_0^T D v\\
\|u^T J_0^T B v\| \le c' r(\|u^2\| + \|v\|^2) & & \|v^T J_0^T C u\| \le c' r(\|u^2\| +
\|v\|^2). 
\end{eqnarray*}
Hence, because $f \ge 1$ always and $f\, v^T = \sqrt f (\sqrt f
v^T)$,  \begin{eqnarray*}
\| u^T J_0^T Bv\|  + \|f\, v^T J_0^T C u\| & \le & 
c' r (\|u^2\| + \|v\|^2) + \sqrt{f} c' r(\|u^2\| + f\|v\|^2)\\
& \le & c' r (1+ \sqrt{f}) (\|u^2\| + f \|v\|^2)\\
& \le & c c'  r (1 + \sqrt{f}) (u^T J_0^T A u + f\, v^T J_0^T D v).
\end{eqnarray*}
Therefore, $\tau_\la$ will tame $J_p$ provided that
$ c c'  r (1 + \sqrt{f} ) < 1$.  Hence  it suffices that 
$$
f = f_\la (r) \le \frac 1{\al r^2}, \;\;\;\mbox{where }\;\;\al = (2 c c')^2
$$
when $r \le r_0$, and $f = 1$ for $r \ge r_0$.  
But for functions of this
type, the integral $$
\int_{D^2} f_\la dx\wedge dy
$$
 can be arbitarily large. The result follows.\QED

The proof in the general case is similar, but it uses a more complicated
normal form for $\tau_0$ near $Z$.  Think of a neighborhood $\Nn(Z)$ of $Z$
as  the unit disc bundle of some complex line bundle over $Z$, let $r$ be
the radial coordinate in this line bundle and choose a connection form $\al$
for the associated circle bundle so that $d\al = -\ka \pi^*(\si_Z)$.
(Here $\ka  = Z\cdot Z > 0$ if we suppose $\si_Z$ normalized to have integral
$1$ over $Z$.)  Then, by the symplectic neighborhood theorem, we may
identify $\tau_0$ with the form $$
\pi^*(\si_Z) + d(r^2\al) = (1 - \ka r^2) \pi^*\si_Z + 2r dr\wedge \al
$$
in some  neighborhood $r \le r_0$.  

We will take $\tau_\la$ to have the form
$$
\pi^*(\si_Z) + d(r^2\al)  - d(f_\la(r)\,\al) = 
(1 - \ka r^2 + \ka\,f_\la) \pi^*\si_Z + (2r - f_\la') dr\wedge \al
$$
for a suitable nonincreasing function $f_\la$ with support in $r\le r_0$.
Thus $f_\la(0) > 0$ and we assume that $f_\la$ is constant very near 
$0$.  The form $- d(f_\la(r) \al)$ represents the
positive multiple $f_\la(0)/\ka$ of $a_Z$.
Hence we need to see that we can choose
$f_\la$ so that $f_\la(0)$ is arbitrarily large and so that $\tau_\la$
tames $J$.  

At each point $p$, split $T_p\Nn(Z)$ into a direct sum $E_H \oplus
E_F$ where the horizontal space $E_H$ is in the kernel of both $dr$ and
$\al$ and where $E_F$ is tangent to the fiber.  Then these subspaces are 
orthogonal with respect to 
$\tau_\la$, for all $\la$ and we may choose bases in these spaces so that
for $(u,v) \in E_H \oplus E_F$,
$$
\tau_0((u,v), (u',v')) =  u^T J_0^T u' + v^T J_0^T v'
$$
as before.  Then
$
\tau_\la((u,v), (u',v')) = a u^T J_0^T u' +
 b v^T J_0^T v',
$
where
$$
a =  1 +  \frac {\ka \,f_\la}{1 - \ka r^2}\;\; \ge\;\; 1,\quad
b =  1 - \frac{f_\la'}{2r} \;\;\ge\;\; 1.
$$
Arguing as above, we find that $\tau_\la$ will tame $J_p$ provided that
 $c \, r (\sqrt a + \sqrt b) < 1$, where  $c$ is a constant
depending only on $J$. 
Such an inequality is satisfied if $ - f_\la' \le c'/r$ for a
suitable constant $c'$ since then $f_\la = const - c' \log(r)$.  It follows that
we can choose $f_\la$ so that $J$ is tamed and $f_\la(0)$ is arbitrarily
large.  Further details are left to the reader.\QED

\begin{cor}\label{cor:infl}  If $J\in \Aa_{\mu}$
admits a curve $Z$ in class $pA + q F$, where $p > 0, q\ge 0,$ then $J\in
\Aa_{\mu'}$ for every $\mu' > q/p$. \end{cor}
\proof{} Choose a form $\tau_0 \in \Ss_\mu$ that tames $J$ and let
$$
\si_\la = \frac 1{1 + \la  p}\;\tau_\la,
$$
where $\tau_\la$ is the family of forms taming $J$ constructed as above.
Then 
$$
\si_\la(F) = 1, \quad \si_\la (A) =  \frac{\mu + \la q}{1 + \la p} = \frac qp +
\eps_\la $$
for arbitrarily small $\eps_\la > 0$.  The result is now immediate.\QED

\section{Finding embedded holomorphic curves}\label{s:emb}

Our aim in this section is to prove Lemma~\ref{le:g0} and~\ref{le:g}.
To do this we have to show that appropriate homology classes $B$ have
embedded $J$-holomorphic representatives for {\it each} given $J$.
Here is the basic method that we will use.

Given a homology class $B$, let $\Gr(B)$ denote the Gromov invariant of $B$
as defined by Taubes.  This counts the  number\footnote
{
This statement holds provided that the class $B$ has no representatives by
multiply covered tori of zero self-intersection. 
It is easy to check that this hypothesis holds in all cases considered here.
Observe also that each curve is equipped with a sign, and that one takes the
algebraic sum.} of embedded $J$-holomorphic curves in class $B$ through
$k(B)$ generic points for a fixed generic $J$, where  $$
k(B) =  \frac 12 (c_1(B) + B^2).
$$
Therefore, whenever $\Gr(B) \ne 0$ we know that there have to be
embedded $J$-holomorphic curves in class $B$ for generic $J$.  
Hence, by Gromov compactness, there is for each $J$ {\it some}
$J$-holomorphic cusp-curve (or stable map)\footnote
{
Here we use the word ``curve" to denote the image of a connected smooth
Riemann surface under a $J$-holomorphic map, and reserve the word
``cusp-curve" or stable map for a  connected union of more than one
$J$-holomorphic curve.
} 
that represents $B$ and goes through an arbitrary
set of $k(B)$ points, and our task is to show that at least one of these
representatives is an embedded curve.  We do this by showing that there are
not enough degenerate representatives to go through $k(B)$ generic points.

One very important point is that, by the work of
Hofer--Lizan--Sikorav~\cite{HLS}, embedded curves in class $B$ are
regular\footnote {
A $J$-curve $u:\Si \to M$ is said to be regular if the linearization $Du$ of the
generalized Cauchy--Riemann operator  at $u$ is surjective.  If this is the
case for all $J$-curves in some moduli space, then this moduli space is a
manifold of the ``correct" dimension, i.e. its (real) dimension equals the index 
$2k(B)$ of $Du$.   Thus the moduli space of $B$-curves with $k(B)$ marked
points has real dimension $4k(B)$, the same as the dimension of $M^{k(B)}$. 
To a first approximation, the Gromov invariant is just the degree of the
evaluation map from this pointed moduli space to $M^{k(B)}$. 
See~\cite{MS1,Mg} 
for further information.} whenever $c_1(B) > 0$.  We
shall need the analog of this result for non-embedded curves, which has been
worked out  by Ivashkovich--Shevchishin~\cite{IS}.   The relevant condition 
is $c_1(B) > g(B)$,
where $g(B)$ is the genus of an embedded representative of $B$.

This method works only because ruled surfaces have many nonzero Gromov
invariants.  Indeed, it was proved by Li-Liu~\cite{LL}  that, if $M = \Si\times
S^2$ where $g > 0$ and $B = pA + qF$, then
$$
\Gr(B) = (p+1)^g, \quad \mbox{provided that }\; k(B) 
\ge 0.
$$
In particular, $\Gr(B) \ne 0$ provided that $q \ge g-1$.  When $g = 0$,
$\Gr(B) = 1$ for all classes $B$ with $p,q \ge 0$ and $ p+q > 0$.
\MS

The following result is well known.  We sketch the proof because most
 references consider only the case of generic $J$.

\begin{lemma}\label{le:fib}  For each $J\in \Aa_\mu$, $M = \Si\times S^2$
admits a smooth foliation by $J$-holomorphic spheres in class $F$.
\end{lemma}
\proof{}  The above remarks imply that  $\Gr(F) \ne 0$ for all $g$.  Hence
the fiber class always has some $J$-holomorphic representative.  When
$g > 0$ any such representative has to be a curve, rather than a cusp-curve,
since $F$ is a generator of the spherical part of $H_2(M)$.  Hence the moduli
space of  $J$-holomorphic curves in class $F$ is compact and has real
dimension $2 k(F) = 2$.  Since $F\cdot F = 0$, the curves have to be disjoint
by positivity of intersections, and they are embedded by the adjunction
formula.   The proof that they form the leaves of a smooth foliation is a little
more subtle and may be found in~\cite{LM1} Prop. 4.12.  See also~\cite{HLS}.

When $g=0$ we must rule out the possibility that $F$ can be represented by
a cusp-curve.  If so, $F$ would decompose as a sum $B_1 + \dots + B_\ell$,
where each $B_i = p_i A + q_i F$ has a spherical $J$-holomorphic
representative and $\ell > 1$.   Thus $\om_\mu(B_i) = \mu p_i + q_i > 0$ for
all $i$.  Since also $\sum_i \om(B_i) = 1$ and $\ell > 1$,  $\om(B_i)< 	1$ for
all $i$.  Therefore, since $\mu \ge 1$, for each $i$ $p_i$ and $q_i$ must both
be nonzero and have opposite signs.  However $B_i\cdot B_j < 0$ if 
$p_i, p_j$ have the same sign.  Therefore, by positivity of 
intersections, $\ell=2$ and $p_1, p_2$ have opposite signs.  Further 
the classes $B_1, B_2$ and hence $F$ have unique representatives, 
contradicting the fact that there is an $F$ cusp-curve through every point.
\QED

\NI
{\bf Proof of Lemma~\ref{le:g0}}

We have to show that when $g=0$ and $J \in \Aa_{\mu, \ell}$ the class $B
= A + \ell F $ has a $J$-holomorphic representative.  By the above, $\Gr(B) =
1$ and $k(B) = 1 + 2\ell$.  Hence there is at least one $J$-holomorphic
representative of the class $B$ through each set of $1 + 2\ell$ points.
Because $A - \ell F$ is represented,  it follows
from positivity of intersections that
no class of the form $A + qB$ with $q <
\ell, q\ne -\ell$ can be represented.  Therefore, the only cusp-curves that
represent $B$ must be the union of the fixed negative curve $C_J$  in class $A
- \ell F$ with $2 \ell $ $J$-fibers.  But if we choose the $1 + 2\ell$ points
so that none lies on $C_\ell$ and no two lie in the same $J$-fiber then there
is no such cusp-curve through these points.  Hence there has to be an
embedded curve through these points.\QED

If one tries the same argument when $g > 0$ we find that
$c_1(B) > 0$ when $\ell > g-1$ and that $k(B) = 1 + 2\ell -g$.  Hence 
there are enough cusp-curves of the form $C_\ell$ plus $2\ell$ fibers
to go through $k(B)$ generic points.    However, now $\Gr(B) = g+1$ and one
might be able to show that these cusp-curves cannot account for all the
needed elements.  For example, when $g = \ell = 1$, there is only one
cusp-curve through $k(B) = 2$ generic points, and so if one could show this
always has ``multiplicity 1" , there would have to be another embedded
representative of the class $B = A+F$.   To carry such an argument through
would require quite a bit of analysis.  Therefore, for now, we will 
only prove the following partial results.

\begin{lemma}\label{le:curp}\begin{itemize}
\item[(i)]
Every $J \in
\Aa_{\mu,1}^1$ admits an embedded representative in 
either $A + F$ or $2A + 3F$.
\item[(ii)]  Every $J \in
\Aa_{\mu,0}$ admits an embedded representative in 
the class  $B = A + qF$ for some integer $q\le g/2$.
\end{itemize}   \end{lemma} 
\proof{}  We first prove $(i)$, so $g  = 1$.
We shall be considering classes
of the form $A + qF$ and $2A + qF$.  Since $(A + qF)\cdot F = 1$, it follows
from positivity of intersections that any $J$-holomorphic  representative $C$
of this class $A + qF$ must either intersect a $J$-holomorphic fiber
transversally in a single point or must contain this fiber completely.  Thus
each component of $C$ must be embedded.  If $C$ represents the class $2A +
qF$, then $C$ could have double points or critical points of order $1$,
i.e. places
where the derivative $du$ of the parametrizing map $u$ has a simple zero.
(The critical points cannot have higher order because $C\cdot F = 2$.)
Since each such critical point of $C$ reduces its genus by $1$ there can be 
at most $g(2A + qF)$ such points, where 
where $g(B)$ denotes the genus of an embedded representative of $B$.   
Thus it follows from Theorem~2 of~\cite{IS} that the curve $C$
 will belong to a smooth moduli
space of the ``correct" dimension provided that $c_1(2A + qF) > g(2A + qF) $.

Now, when $g = 1$,
$$ 
\begin{array}{lll}
k(A + qF) = 2q, & c_1(A + qF) =  2q, & g(A + qF) = 1,\\
k(2A + qF) = 3q, & c_1(2A + qF) = 2q, & g(2A + qF) = 1+q.
\end{array} $$
Therefore $c_1(2A + qF)  > g(2A + qF)$ provided that $q > 1$.
In this case $k(2A + qF) > 0$ so that $\Gr(2A + qF) > 0$.

Because the classes $F$ and $A - F$ have $J$-representatives by
hypothesis,  the only other classes $pA + qF$ that can be represented have
$p \ge 0$ and $q \ge p$.
Thus, the possible representatives for
 the class $2A + 3 F$ have  one
of the following types.\MS

\NI
(a)  an embedded representative of   $2A + 3 F$;
\SmS

\NI
(b)  a representative of   $2A + 3 F$ with some critical or double points;
\SmS

\NI
(c)  the union of two embedded curves in class $A + q_i F$, where $q_i\ge
1$,  with $3 - q_1 - q_2$ fibers;\SmS

\NI
(d)  the union of the $A - F$ curve, together with an embedded $A + qF$
curve,  with $q\ge 1$, and $4 - q$ fibers;\SmS

\NI
(e)  a $2$-fold copy of the $A -  F$ curve together with $5$
fibers.
\MS

 We claim that 
if (a) does not hold, there is an embedded curve in class $A + F$. 
Suppose not.  Then there can be no cusp-curves of type (c) since one of the
$q_i$ would have to be $1$.   Therefore it remains to show that there
are not enough  cusp-curves of types (b), (d) or  (e) to go through
$k(2A+3F) = 9$  generic points. 

Let us consider these types in turn.  Since
$c_1(B) > g(B)$ when $B = 2A + mF$ any curve of type (b) is regular and
so belongs to a moduli space of (complex) dimension strictly smaller than
$9$: indeed, it is easy to see from~\cite{IS} that the difference is
precisely the number of critical points.  

Now consider cusp-curves of type (d).   In this case,  the $A + qF$ curve
would be regular since $q > 1$.  Hence this
cusp-curve could go through at most
$$
k(A + qF) + 4 - q \;\le\; k(A + 4F) = 8
$$
points.  
Finally, curves of type (e) go through at most $5$ generic points. This proves
$(i)$.

The proof of $(ii)$ is easier.  Let $m = k$  where $g = 2k+1$ or $2k$.
Then $k(A + mF) \ge 0$. Therefore the class $A + mF$ has to have some
$J$-holomorphic representative.  If $J\in \Aa_{\mu, 0}$, this must be the
union  of an embedded curve in class $A + kF$ with $m-k$ fibers for some $k
\le m$ since there are no negative $J$ curves.\QED

\begin{cor}  Lemma~\ref{le:g} holds.
\end{cor}
\proof{}  Apply Proposition~\ref{prop:key}.\QED

\begin{rmk}\label{rmk:main}\rm
The argument in Lemma~\ref{le:curp} can be generalized to the case $g > 1$.
The results are not very sharp because  the
 requirement that $c_1(B) > g(B)$,
which is used  to ensure that the curve is embedded, forces us to consider
only classes $B = pA + qF$ with $p \le 2$.   However, it is quite possible that
we might be able to find an embedded $J$-curve under weaker hypotheses. 
For example, if $\ell \ge g$, any $J$-curve in class $B$ would have to have
$q/p \ge \ell > g-1$ which is enough to imply that $c_1(B) > 0$.   Hence all
the $J$-curves in class $B$ would be regular.  It is then not hard to adapt the
above argument to show that there is a nonempty moduli space of curves in
some such class $B$ with $q/p$ arbitrarily close to $\ell$.  It is possible that
one could show that a generic element of this moduli space must be
embedded. If so, then we would be able to prove the analog of part $(i)$ of
Theorem~\ref{thm:main} when $\ell > g$.
\end{rmk}

\section{The groups $G_\mu^g$}\label{s:dd}

We first study the topology of the groups $\Dd^g$ and $\Dd_0^g$, and 
then prove part $(ii)$ of Proposition~\ref{prop:infty}.  

\subsection{The groups $\Dd^g$ and $\Dd_0^g$.}\label{ss:dd0}

Consider the group $\Dd^g$ of
fiber preserving diffeomorphisms of $\Si\times S^2$, and let $\io:\Dd^g\to
\Diff(\Si\times S^2)$ be the inclusion. 

\begin{lemma}\label{le:pi0}  The map $\io_*: \pi_0(\Dd^g) \to 
\pi_0(\Diff(\Si\times S^2))$ is injective.
\end{lemma}
\proof{}  Clearly,
$$
\Dd^g \simeq \Diff(\Si)\times \Map(\Si, \SO(3)).
$$
Since $\pi_2(\SO(3)) = 0$, any map $\Si\to \SO(3)$ that induces the zero map
on $\pi_1$ is homotopic to a constant map.  Hence
$$
\pi_0(\Dd^g) = \pi_0(\Diff(\Si)) \times \Hom(\pi_1(\Si), \SO(3)) =
\pi_0(\Diff(\Si)) \times (\Z/2\Z)^{2g}. 
$$
Next observe that there are maps
$$
\Diff(\Si\times S^2) \to \Map(\Si\times S^2, S^2) = \Map (\Si,\Map_1(S^2,
S^2))\to  \Map(\Si, \SO(3)),
$$
where the first map takes a diffeomorphism $\phi$ to its composite with the
projection onto $S^2$, and the second exists because there is a projection
of the space $\Map_1(S^2, S^2)$ of degree $1$ self-maps of the sphere
onto $\SO(3)$ (see~\cite{A}.)  It follows that the part of $\pi_0(\Diff(\Si))$
coming from $\Map(\Si, \SO(3))$ injects into $\pi_0(\Diff(\Si\times S^2))$.
The rest of  $\pi_0(\Diff(\Si))$ can be detected by its action on the set of
conjugacy classes in the group $\pi_1(\Si\times S^2)$, and hence also
injects.\QED

The following corollary was used in~\S\ref{ss:bi}.

\begin{cor}\label{cor:con} $\Dd^g\cap \Diff_0(\Si\times S^2) = \Dd_0^g$.
\end{cor}

Now consider the identity component $\Dd_0^g$.  Clearly
$$
\begin{array}{lll}
\Dd_0^0 \simeq \SO(3) \times \Map_0(T^2, SO(3)), & \mbox{ when } &\Si =
S^2,\smallskip\\ 
\Dd_0^1 \simeq T^2 \times \Map_0(T^2, SO(3)), & \mbox{ when } &\Si =
T^2,\smallskip\\ 
\Dd_0^g \simeq \Map_0(\Si, SO(3)), & \mbox{ otherwise},&
\end{array}
$$
where $\Map_0$ denotes the component containing the constant maps.
Consider the evaluation map
$\Map_0(\Si, \SO(3)) \to
\SO(3)$.
 Because the elements of $\SO(3)$ lift to the
constant maps, $\Map_0$ is a product: $$
\Map_0(\Si, \SO(3)) \cong \SO(3)\times \Map_{*}(\Si, \SO(3)),
$$
where  $\Map_{*}$   is the connected space  of all based
null-homotopic maps. 
Moreover, since every such map $\Si\to \SO(3)$ has a
unique lift to $S^3$,
$$
 \Map_{*}(\Si, \SO(3)) =  \Map_{*}(\Si, S^3).
$$
We write $\Om(X) = \Map_*(S^1, X)$ for the identity component of the based
loop space of $X$, and $\Om^2(X) = \Map_*(S^2, X)$ for the identity
component of the double loop space.

\begin{lemma}  $\Map_{*}(\Si, S^3)$ is homotopy equivalent to
the product of $\Om^2(S^3)$ with $2g$ copies of
$\Om(S^3)$.
\end{lemma}
\proof{} Choose loops $\ga_1,\dots, \ga_{2g}$ in $\Si$ that represent 
 a standard basis for
$H_1(\Si)$ and choose corresponding projections $pr_j:\Si\to S^1$  such that
$pr_j$ is injective on $\ga_j$ and maps each $\ga_i, i\ne j,$
to the base point $id \in S^3$.  Each element $h\in \Map_{0,*}(\Si, S^3)$
determines an element  $(f_i) = (h|_{\ga_i}\circ pr_j)$ of
the $2g$-fold product $(\Om(S^3))^{2g}$ by restriction to the loops $\ga_i$. 
Then, if $\cdot$ denotes the product of maps coming from the group
structure on $S^3$, the map
$$
h' = h\cdot f_1^{-1}\cdot\dots\cdot f_{2g}^{-1}
$$
takes each loop $\ga_i$ to the base point in $S^3$ and so can be considered
as an element of $\Om^2(S^3)$.  It is not hard to check that the map
$h \mapsto (f_i, h')$ is a homotopy equivalence.\QED

\begin{cor}  The vector space $\pi_i(\Dd_0^g)\otimes \Q$ has
dimension $1$ when $i = 0,1, 3$  except in the cases $i = g =
1$ when the dimension is $3$ and $g = 0, i = 3$ when the
dimension is $2$.  It has dimension  $2g$ when $i = 2$, and is zero otherwise.
\end{cor}

\begin{cor}\label{cor:d0}  When $\mu > 1$ the map $G_\mu^0 \to \Dd_0^0$
induces a surjection  on rational homotopy.
\end{cor}
\proof{}  This follows immediately from Proposition~\ref{prop:am3}, since
the elements in $\pi_3(\Dd_0^g)\otimes \Q$  lift to $G_\mu^0$ for all $\mu$.
\QED

\begin{cor}\label{cor:surj}  For all $g > 0$
 the map $G_\mu^g\to \Dd_0^g$
induces a surjection on $\pi_2\otimes \Q$ for all $\mu > 0$. \end{cor}
\proof{}  We just have to exhibit explicit
representatives in $G_\mu^g$ of the elements in $\pi_2$ for any $\mu> 0$. 
For definiteness, consider the case $g=1$. Then the generator of
$\pi_2(\Dd_0^1)$ corresponding to the projection of $T^2$ onto its first factor
is  $$
S^2\times T^2 \to \SO(3): (z,s,t)\mapsto R_{z,s}
$$
where $R_{z,s}$ is the rotation  of $S^2$ by angle $s$ about the axis through
$z$.  Consider the corresponding family of diffeomorphisms
$$
\rho_z: T^2\times S^2 \to T^2\times S^2:\qquad (s,t,w)\mapsto
(s,t,R_{z,s}w),\quad z\in S^2.
$$ 
Then $\rho_z^*(\om_\mu) = \om_\mu + ds\wedge \al_s$ where $\al_s $
is an exact $1$-form on $S^2$ or each $s$.  Thus
$\rho_z^*(\om_\mu) $ can be isotoped to
$\om_\mu$  via the symplectic forms $\om_\mu + \nu ds\wedge \al_s$,
where $\nu\in [0,1]$.  It follows, using Moser's argument, that the family
$\rho_z$ can be homotoped into $G_\mu^1$ for any $\mu$.\QED

\begin{rmk}\rm\label{rmk:homog}
The whole analysis of the groups $G_\mu^0$ is based on  Gromov's
result that when the two spheres have equal size (i.e. when $\mu = 1$) 
$G_1^0$ is homotopy equivalent to the stabilizer $\Aut(J_{split}) =
\SO(3)\times \SO(3)$ of the product almost complex structure.  
More precisely he showed that the (contractible) space $\Jj_1$ is homotopy 
equivalent to the coset space $G_1^0/\Aut(J_{split})$.  It is shown
in~\cite{AM} that this continues to hold for the other strata $\Jj_{\mu,\ell}$
in $\Jj_\mu$, i.e. each stratum $\Jj_{\mu,\ell}$  is homotopy equivalent to
some homogeneous space $G_\mu^0/ K_\ell$, where the Lie group $K_\ell$ is
the stabilizer of an integrable element in $\Jj_{\mu,\ell}$.   The basic reason
for this is that  the topological type of the space of all $J$-holomorphic curves
 does not change as  $J$ varies  in a
stratum $\Jj_{\mu, \ell}$.

It is conceivable that when $g > 0$
a similar result continues to hold for some of the $\Jj_{\mu,\ell}$, for
example if $\mu$ and $\ell$ are sufficiently large.  However, it does not hold
when $\ell = 0$, $\mu \le 1$.  In this case, because $\Jj_{\mu, 0} = \Jj_\mu$
is contractible, such a statement would imply that $G_\mu$ has the
homotopy type of the stabilizer of some almost complex structure $J$ and so
would be homotopy equivalent to a Lie group.   This is not true because 
$\pi_2(G_\mu)\otimes \Q \ne 0$.  Note also that in this case 
the topological type of the space of $J$-holomorphic curves
may change as  $J$ varies  in $\Jj_{\mu, 0}$.  For example, when $g = 1$,
although the Gromov invariant of the class $A = [T^2\times pt]$ is $2$, there
can be arbitrarily many pairs of $J$-holomorphic tori in class $A$ that
cancel out in the algebraic sum: see the
discussion at the end of  [M3]\S1.
\end{rmk}

\subsection{Relations between $G_\mu $ and $\Dd_0^g$}\label{ss:dd}

It remains to prove part $(ii)$  of Proposition~\ref{prop:infty}
that claims that the map $G_\mu^g\to \Dd_0^g$ induces a surjection on 
rational homotopy when $g> 0, \mu > 0$.
It follows easily from the results in \S\ref{ss:dd0} (see particularly
Corollary~\ref{cor:surj}) that this holds in all dimensions except possibly 
for $1$.   Hence it suffices to show:

\begin{prop}  The map $\pi_1(G_\mu^g)\otimes
\Q \to \pi_1(\Dd_0^g)\otimes
\Q$ is surjective for all $g> 0, \mu > 0$.
\end{prop}
\proof{}    Recall that $\pi_1(\Dd_0^g)\otimes \Q$ has dimension $1$ except
when $g = 1$ when there are two extra dimensions coming from 
translations of the base.   Since these translations lie in $G_\mu^1$ for all
$\mu$, we only need to show that the generator of 
$\pi_1(\Dd_0^g)$ coming from $\pi_1(\Om^2(S^3))$ lifts to the $G_\mu^g$.

As mentioned in \S1.1 when $g=0$ this
generator   can be  explicitly represented by
the circle action
$$
\phi_t:  S^2\times S^2 \to S^2\times S^2:\quad
(z,w)\mapsto (z, R_{z,t}(w)), \qquad t\in [0,1],
$$
where $R_{z,t}$ is the rotation of the fiber sphere $S^2$ through the angle
$2\pi t$ about the axis through the point $z\in S^2$.   Hence the
corresponding generator of $\pi_1(\Dd_0^g)$  can be  represented by 
the circle
$$
\phi_t^g:  \Si\times S^2 \to \Si\times S^2:\quad
(z,w)\mapsto (z, R_{\rho(z),t}(w)),\quad t\in [0,1]
$$
where $\rho:\Si\to S^2$ has degree $1$.  Therefore, we just have to see that
 $\{\phi_t^g\}$ can be homotoped into $G_\mu^g$ for all
$g> 0, \mu > 0$, i.e. that the orbits of its action on $\Ss_\mu^g$
are contractible. 
Equivalently, we need to see  that its orbits 
contract in $\Aa_\mu^g$ for all $ g > 0, \mu > 0$.   

In fact, we shall see that the circle action on $\Aa_\mu^g$ has a fixed point
$J_0^g\in \Aa_{\mu,1}^g$ when $\mu > 1$, so that its orbits obviously
contract in this case.  If we choose $J\in \Aa_{\mu,0}^g$ close to $J_0^g$ its
orbit lies in some normal slice $D$ to the stratum $\Aa_{\mu,1}^g$ at $J_0^g$.
Since $D$ has dimension at least $4$ when $g > 0$, the orbits are contractible
in $\Aa_{\mu,0}^g$ also.  This completes the proof when $g = 1,$ since
Lemma~\ref{le:g} says that
the sets $\Aa_{\mu,0}^1$ are the same for all $\mu > 0$.    However, 
the latter statement may not hold in higher genus and so
the proof  requires more work.  The idea is to pull a contracting
homotopy back from $\Aa_{\mu, 0}^1$ to $\Aa_{\mu, 0}^g$.  

If $g' > g$ an  arbitrary  almost complex structure $J$ on
$\Si_g\times S^2$ does not of course lift to $\Si_{g'}\times S^2$.  However,
suppose that $J$  splits as a product 
on some open neighborhood  $U\times S^2$ of some fixed fiber $F_* = *\times
S^2$ and choose a smooth map $\si:\Si_{g'}\to \Si_g$ that is bijective except
over  the set $\si^{-1}(*)$.  Then, if $\overline V$ is a small closed disc
centered on $*$ that is contained in  $U$,  there is an almost
complex structure $J'$ on $\Si_{g'}\times S^2$ that equals the pullback of $J$
via $\si\times id$ on $\si^{-1}(T^2 - \overline V)$ and  is a product on
$\si^{-1}(U)\times S^2$.  Note that $J'$ is tamed by a symplectic form of the
type
$$
\om' = (\si\times id)^*(\om) + pr^*(\tau)
$$
where $\tau$ is a nonnegative $2$-form on $\Si_{g'}$ with support near
$\si^{-1}(*)$ and $pr$ is the projection.   Thus $J'\in \Aa_{\mu'}^{g'}$ for some
$\mu' > \mu$.  The second condition below is introduced to give control over
$\mu'$.

We will say that $J\in \Aa_{\mu, 0}^g$ is {\it normalized} if satisfies the
following conditions:\MS

\NI
$(i)$  it splits as a product near $F_*$,\smallskip

\NI
$(ii)$ there is a $J$-holomorphic curve $Z$ in class $A$ that is flat near
$F_*$, i.e. it coincides with some flat section $\Si_g\times pt$ near $F_*$.\MS

As explained above, the first  condition allows $J$ to be pulled back to
an element $J'$ of $\Aa_{\mu'}^{g'}$ for some $\mu' > \mu, g' > g$.  Moreover,
$J'$ is necessarily in $\Aa_{\mu', 0}^{g'}$, since it clearly  admits a curve $Z'$
(the pullback of $Z$) in class $A$.  Inflating along $Z'$ we find  that
 $J' \in \Aa_{\eps, 0}^{g'}$ for all $\eps > 0$.   

We construct below,  for some $\mu > 1$, a $2$-dimensional  family $J_w,
w\in D,$ of  elements in $\Aa_{\mu}^0$ that intersect the stratum 
$\Aa_{\mu,1}^0$ transversally at $w = 0$ and are normalized for $w\ne 0$
via the curves $C_w$.  This normalization is uniform in the sense that 
both conditions $(i)$ and $(ii)$ are satisfied over
 some fixed neighborhood of the fiber $F_*$.  
Thus the central element $J_0$ also satisfies condition $(i)$ above, and as
$w\to 0$ the curves $C_w$ converge to the union of the $J_0$-holomorphic
curve $C_0$ in class $A-F$ with a fixed fiber $F_0$ distinct from $F_*$
where $C_0$ is  flat near $F_*$.    Moreover,
$\phi_t^0$ acts on the $J_w$ via   $$
(\phi_t^0)_* (J_w) = J_{e^{it} w}.
$$

This construction permits us to complete the proof as follows. 
We can define the pullbacks $(C_w^1, J_w^1)$  of the pairs $(C_w, J_w)$ 
to $T^2\times S^2$ so
that  the $J_w^1$  all lie in some set
$\Aa_{\mu'}^1$ and the circle acts as above.  By
construction the  $2$-disc $\{J_w^1: w\in D\}$ is normal to the stratum
$\Aa_{\mu,1}^1$, and so sits inside a transverse slice of dimension $4$
consisting of elements $J_{v}^1$ where $v = (w,w')\in D'' =  D\times
D'$.  By the
remark at the end of the proof of Lemma~\ref{le:str} we
can assume that $J_{(w,w')}^1 = J_{w,0}^1$ near $F_*$, so that the
$J_{v}^1$
satisfy condition $(i)$ above .  Moreover, assuming as we may that $J_0$ is a
generic complex structure on the torus $C_0$, we can extend the family of
$J_{w}^1$-holomorphic curves $C_w^1$ to all small $v\in D''$ by the process
of  gluing as in~\cite{M3}.\footnote
{Note that we
need $g=1$ here in order for the problem to have nonnegative index: when
$g > 1$ generic elements in $\Aa_{\mu}^g$ do not admit curves in class $A$.}
Here we are making the curves $C_{v}^1$  by gluing $C_0$ to the fiber
$F_0$.  Therefore, because  the $C_{v}^1$ converge $C^\infty$ to $C_0\cup
F_0$  away from the point $C_0\cap F_0$, they  will be approximately flat
 near $F_*$ and will become flatter as $|w'|\to 0$.    

Now suppose that we want to show that the orbits of $\{\phi_t^g\}$ contract
in $\Aa_\mu^g$ for some particular $g$ and $\mu > 0$.
We can pull the family $(C_{v}^1, J_{v}^1)$ back to a family 
 $(C_{v}^g, J_{v}^g)$ where $J_{v}^g\in \Aa_{\mu'}$ for all ${v}\in D$.
The curve $C_{v}^g$ will now be only 
 approximately $J_{{v}}^g$-holomorphic, but  its deviation
from being holomorphic will tend to $0$ as $|w'|\to 0$.
This means that  inflating along
$C_{v}^g$ will give rise to symplectic forms $\om_\la \in \Ss_\la$ that
tame  $J_{v}^g$ for some restricted interval $\la_0 \le \la \le \mu'$,
where we can make $\la_0$ as close to $0$ as we want
by making $|w'|$ sufficiently close to $0$.   Define
$$
W_\de = \{v \in D'' - \{0\}: |w'|\le \de\}.
$$
The above discussion shows that for any $\mu > 0$ there is $\de > 0$ so
that the pullbacks $J_{v}^g$ of  the elements $J_{v}^1, v\in W_\de$ 
lie in $\Aa_\mu^g$.   By construction, this set  $\{J_{v}^g: v\in W_\de\}$
contains an orbit of $\{\phi_t^g\}$.   
 Since  $W_\de$ is simply connected, this orbit contracts in 
 $\{J_{v}^g: v\in W_\de\}$ and hence in  $\Aa_\mu^g$.

It remains to construct the family $(J_w, C_w)$ on
$S^2\times S^2$  for $w \in D\subset \C$. To do this, take two  copies of
$\C\times \CP^1$  and glue them via the map
$$
\left((\C - \{0\})\times \CP^1\right)_1 \to \left((\C - \{0\})\times
\CP^1\right)_2:\quad (z, [v_0: v_1])_1\mapsto (\frac 1z, [z^2 v_0 + w zv_1:
v_1])_2.   $$
Identify the resulting complex manifold $M_w$  with 
$(S^2\times S^2, J_w)$ so as to preserve the product structure near $F_* =
(\{0\}\times \CP^1)_2$.   For each $w\ne 0$ and $a\in
\C$, there is a holomorphic section $Z_{a,w}$ of $M_w$ containing the points
$$
(z, [1: a - \frac zw])_1, \quad  (u, [wa: au - \frac 1w])_2,
$$
where $u = 1/z$.
As $a$ varies (and $w$ is fixed) these are disjoint.  Hence they must lie in
class $A$.   Moreover $C_w = Z_{0,w}$ contains the points $(u, [0:1])_2$ and so
its image in $S^2\times S^2$ is flat near $F_*$. 
Finally, if we make  $\phi_t^0$  act on each copy of
$\C\times \CP^1$ via $(z, [v_0: v_1])_j\mapsto (z, [v_0: e^{-it} v_1])_j$ and 
then it 
takes $M_w$ to $M_{e^{it}w}$.
 The other properties are obvious.
\QED

\begin{cor}\label{cor:const}  For all $i$ and $g > 0$, the kernel of the map 
$$
\pi_i(\Diff_0(\Si\times S^2))\otimes \Q
\to \pi_i(\Ss_\mu)\otimes \Q
$$ 
is independent of $\Q$.
\end{cor}
\proof{}  The kernel of this map is the image of $\pi_i(G_\mu^g)$.  The result
follows from the previous lemma because the inclusion $G_\mu^g \to
\Diff_0(M)$ factors through $\Dd_0$.\QED

\end{document}